\documentclass[11pt,a4paper]{article}

\usepackage[margin=25mm]{geometry}
\usepackage[T1]{fontenc}
\usepackage[utf8]{inputenc}
\usepackage{lmodern}
\usepackage{amssymb,amsthm,mathtools,mathrsfs}
\usepackage{microtype}
\newcommand{\T}{\mathbb{T}}
\newcommand{\Es}{E^{\mathrm{s}}}
\newcommand{\dd}{\,\mathrm{d}}
\usepackage{enumitem}
\usepackage{hyperref}
\hypersetup{hidelinks}
\newtheorem{theorem}{Theorem}[section]
\newtheorem{lemma}[theorem]{Lemma}
\newtheorem{proposition}[theorem]{Proposition}
\newtheorem{corollary}[theorem]{Corollary}
\theoremstyle{definition}

\theoremstyle{remark}
\newtheorem{remark}[theorem]{Remark}

\numberwithin{equation}{section}
\newenvironment{keywords}{\par\smallskip\noindent\textbf{Keywords:}\ }{\par\smallskip}
\newenvironment{MSCcodes}{\par\noindent\textbf{2020 MSC:}\ }{\par\medskip}

\newcommand{\R}{\mathbb R}

\newcommand{\norm}[1]{\left\|#1\right\|}

\newcommand{\ip}[2]{\left\langle #1,#2\right\rangle}
\allowdisplaybreaks[2]
\setlist[enumerate]{label=\textup{(\roman*)},leftmargin=2.1em,itemsep=2pt,topsep=3pt}

\title{Realization of expanding dynamics by a finite horizon parabolic law response}
\author{Masood Shahbaz Shoukat\thanks{Guizhou University, Huaxi Park, Guiyang 550000, Guizhou, China.}
\and Luca Di Persio\thanks{University of Verona, Strada Le Grazie, 37134 Verona, Italy.}
\and Wensheng Jia\thanks{Guizhou University, Huaxi Park, Guiyang 550000, Guizhou, China. Corresponding author (\href{mailto:wsjia@gzu.edu.cn}{\texttt{wsjia@gzu.edu.cn}}).}}
\date{}

\begin{document}
\maketitle
\begin{abstract}
We study the discrete dynamical system generated by iterating an episodic response of a population density. In each episode, the running potential is frozen at the current density, agents solve a quadratic control problem, and the terminal density is passed to the next episode. For each prescribed positive horizon $\tau$ and each $C^k$ map $F$ of finitely many Fourier moments, with $1\le k<\infty$, we construct a finite mode interaction realizing $F$ at order $\varepsilon$. In scaled coordinates, the exact parabolic response satisfies $\|T^{(d)}_{\varepsilon,\tau}-(F,S(\tau)|_{\Es_d})\|_{C^k}\le C\varepsilon$. The estimate is uniform over a full bounded infinite dimensional Sobolev fibre and needs neither Fourier truncation nor an invariant manifold. For scalar targets with disjoint expanding branches that strictly cover their convex hull, this $C^1$ control and parabolic compactness yield compact invariant sets of the exact response. These sets factor onto the full one sided shift and contain periodic densities of every least period. Nondegenerate periodic orbits continue, and hyperbolicity requires a separate strict multiplier condition. We also prove a derivative cost obstruction for nonaffine microscopic realization, a susceptibility weighted criterion for population energy designs, and the failure of exact scalar closure under nontrivial rank one feedback. The results concern designed frozen law responses, not Nash consistent mean field game equilibria or finite player limits.
\end{abstract}


\section{Introduction}\label{sec:intro}

We study the dynamical system generated by iterating a \emph{frozen law optimal control response}. At the start of each episode, a running potential is evaluated at the current population density and held fixed until the terminal time $\tau$. Agents then solve the resulting quadratic-control problem. This is a Hamilton--Jacobi--Bellman/Fokker--Planck pair, which the Cole--Hopf transformation reduces to linear parabolic equations. The terminal density becomes the input of the next episode. The structure is that of quadratic mean field games \cite{LL,HMC,CD,Gueant,SGU,USG}, where individual optimality is coupled to an evolving population law. Here, however, no equilibrium condition couples the running cost to the evolving law inside an episode. Agents optimize over a finite horizon against a frozen snapshot of the population law, so the response is not a Nash consistent forward--backward equilibrium. The iteration of this response map is the dynamical system studied here.

Write $e_j(x)=\cos(2\pi jx)$ on $\T=\R/\mathbb Z$, let $c_j(m)=2\int_\T m e_j$, and fix a viscosity $\nu>0$ and a horizon $\tau>0$. The heat multiplier and the linear susceptibility in mode $j$ are
\begin{equation}\label{eq:introchi}
 E_j=e^{-\nu(2\pi j)^2\tau},\qquad
 \chi_j=\frac{(1-E_j)^2}{2\nu^2(2\pi j)^2}>0.
\end{equation}
For a prescribed profile $f$, the feedback
\begin{equation}\label{eq:introfb}
 V_{\varepsilon,\tau}[m](z)
 =\varepsilon\frac{E_1 a-f(a)}{\chi_1}\,\cos(2\pi z),
 \qquad a=\frac{c_1(m)}{\varepsilon}\,,
\end{equation}
produces the scalar response $f(a)+O(\varepsilon)$ and damps the remaining modes. This formula illustrates the inverse-design mechanism: the susceptibility $\chi_1$ cancels the discrepancy between the heat response $E_1a$ and the prescribed target $f(a)$. The following full-space estimate makes this quantitative for finitely many selected modes.

Fix $\nu>0$, a prescribed horizon $\tau>0$, $s>3/2$, an integer $d\ge1$, an open set $U\subset\R^d$, a compact set $K\Subset U$, and $F\in C^k(U;\R^d)$ with $1\le k<\infty$. For the finite mode interaction \eqref{eq:multifb}, let $T^{(d)}_{\varepsilon,\tau}$ denote the exact parabolic density response in the scaled coordinates $m=1+\varepsilon\big(\sum_{j\le d}a_je_j+y\big)$, $(a,y)\in K\times\overline B_1(\Es_d)$. Theorem~\ref{thm:multimode} states that
\begin{equation}\label{eq:intro-main}
 \big\|T^{(d)}_{\varepsilon,\tau}-\big(F,\,S(\tau)|_{\Es_d}\big)\big\|_{C^k(K\times\overline B_1(\Es_d))}\le C\varepsilon .
\end{equation}
Here $S(t)=e^{\nu t\partial_{zz}}$, and $\Es_d$ is the complement of the first $d$ cosine modes in the mean-zero even Sobolev space $H^{s-1}$; the precise spaces and norms are defined in Section~\ref{sec:main}. For $d=1$ and $F=f$ this is the map $T_{\varepsilon,\tau}$ of \eqref{eq:coords}, so Theorem~\ref{thm:scalar} is the scalar case. The map $T^{(d)}_{\varepsilon,\tau}$ is the exact density update, not a Fourier truncated evolution. The bound controls the selected moments \emph{and} the entire complementary Sobolev component, with derivatives through order $k$, uniformly over the full infinite dimensional stable ball; no finite dimensional invariant manifold and no exact closure of the selected moments is assumed.

In the scalar case, suppose that $f$ satisfies the strict expanding-covering hypotheses of Theorem~\ref{thm:symbolic}. The limiting map is $T_0(a,y)=(f(a),S(\tau)y)$ with $\|S(\tau)|_{\Es}\|=e^{-\nu(4\pi)^2\tau}<1$. For sufficiently small $\varepsilon$, the full-space $C^1$ estimate preserves the strict endpoint crossings, unstable expansion, stable trapping and graph-cone inequalities required by Lemma~\ref{lem:cover}. Relative compactness of the image is supplied separately by positive-time parabolic smoothing (Proposition~\ref{prop:analytic}). The covering lemma therefore yields a compact strictly invariant set of the exact response and a continuous surjection onto the full one-sided shift that intertwines the two dynamics; its periodic-closing argument also produces periodic densities of every least period. Thus the approximation estimate is quantitative, while the invariant-set conclusion concerns the exact PDE-generated density map, and no exact autonomous closure of the first moment is required (Proposition~\ref{prop:closure}). This conclusion is established for scalar targets with expanding interval branches; Theorem~\ref{thm:multimode} alone does not assert it for arbitrary multidimensional targets, and a full-shift factor is weaker than a conjugacy.

The construction combines the exact parabolic representation of the frozen HJB--Fokker--Planck response with an explicit diagonal susceptibility obtained by Duhamel differentiation. The designed feedback inverts this susceptibility on the selected modes. Uniform analytic remainder estimates then give $C^k$ control over the full infinite-dimensional Sobolev fibre, while positive-time smoothing supplies the compactness needed for the invariant-set construction.

The horizon is arbitrary but fixed: for each prescribed $\tau>0$ the interaction is tuned with the corresponding susceptibilities $\chi_j(\tau)$, and the conclusions hold for all $0<\varepsilon<\varepsilon_0(\tau)$ with constants $C=C_\tau$. No long horizon limit and no choice $\tau=\tau(\varepsilon)\to\infty$ is used, and no uniformity as $\tau\downarrow0$ is claimed.

Our contributions are of three kinds. First, we prove the realization estimate \eqref{eq:intro-main} (Theorems~\ref{thm:multimode} and~\ref{thm:scalar}): for every finite $k\ge1$, $C^k$ scalar targets, and $C^k$ maps of finitely many Fourier moments, are realized to order $O(\varepsilon)$ in $C^k$ by the exact response, uniformly over the infinite-dimensional stable ball. Second, we obtain dynamical consequences for the exact PDE (Theorems~\ref{thm:symbolic}--\ref{thm:robust} and Corollary~\ref{cor:examples}): strictly expanding coverings yield compact invariant sets of the exact response factoring onto the full $N$-shift, with periodic points of all periods; a nondegenerate target cycle continues uniquely, with exactly one unstable multiplier when the target cycle is expanding; and the full-shift conclusion is robust under small perturbations of the interaction in the natural scaled $C^1$ norm. Third, we prove restrictions on the designs (Theorems~\ref{thm:cost} and~\ref{thm:potential} and Proposition~\ref{prop:closure}): interactions with bounded first and second population derivatives cannot realize a nonaffine scalar map on shrinking tubes, and the second-derivative cost of order $\varepsilon^{-1}$ is sharp when first derivatives stay bounded; among finite-mode targets, those generated by a population energy are characterized by exactness of a susceptibility-weighted one-form; and for nontrivial rank-one feedback, an explicit coefficient shows that the first moment is not an exactly autonomous factor.
Appendix~\ref{app:long} treats a horizon-independent variant at long horizons.

The construction is an inverse-design result: the coupling is engineered from $f$, and it is not a translation-invariant convolution with fixed Fourier coefficients. For a nonaffine microscopic limit, Theorem~\ref{thm:cost} gives a tradeoff between the first and second population derivatives of the interaction; when the first derivatives stay uniformly bounded, the second derivatives must grow at least at order $\varepsilon^{-1}$, and the explicit rank-one design attains this order whenever $f''\not\equiv0$ on $J$ (see Remark~\ref{rem:cost} for the quantifiers). The estimates are local in amplitude, and depend on the target, the number of modes and the horizon; none is uniform as $\tau\downarrow0$, $d\to\infty$, or as the symbolic alphabet grows. The mean-field terminology refers only to the optimal-control provenance of the response: the Hamilton--Jacobi--Bellman/Fokker--Planck pair is solved with the population cost frozen at the initial law, and no consistency condition is imposed (Appendix~\ref{app:control}). The analysis takes place in $H^r$ with $r>1/2$, so it does not address atomic empirical measures, finite-player limits, attracting physical measures, or topological conjugacy with a Cantor shift. Valid dynamical conclusions rely on strict covering, nondegeneracy, or spectral separation, as appropriate. The distinction between continuation and hyperbolicity is sharp: for $f(a)=-a$ the fixed point is nondegenerate yet has multiplier exactly $-1$ for every $\varepsilon$ (Remark~\ref{rem:nonhyp}).

Prescribing dynamics inside dissipative PDEs has a long history. Pol\'a\v cik \cite{Pol} and Pol\'a\v cik--Rybakowski \cite{PR} embedded finite-dimensional vector fields in scalar parabolic equations, and Fiedler and Pol\'a\v cik \cite{FP} obtained complicated dynamics in scalar reaction--diffusion equations with a spatial-average nonlocal term. Those results concern invariant continuous-time dynamics of specially constructed equations; ours concerns a discrete map generated by a fixed-horizon two-equation response. Self-consistent transfer operators \cite{Gal,CGT,CGT26} provide a framework in which a population update is a nonlinear operator on an infinite-dimensional space; they do not supply the realization formula proved here. Bifurcations in a reduced finite-horizon mean-field model are studied in \cite{Lori}. The covering, symbolic-coding and graph-transform arguments are classical \cite{HPS,ZG,Walters}, and we make no claim that they or the realization of prescribed dynamics in parabolic equations are new in themselves. The exactness test behind Theorem~\ref{thm:potential} is likewise classical \cite{MS,Letcher}; what is specific here is the diagonal, strictly positive finite-horizon susceptibility that sets the weights of the integrability condition.

After this introduction (Section~\ref{sec:intro}), Section~\ref{sec:main} states the model and proves the realization, dynamics and obstruction results. Section~\ref{sec:conclusion} discusses scope and open problems. The appendices contain the functional-analytic facts (Appendix~\ref{app:analytic}), the covering lemma (Appendix~\ref{app:cover}), the long-horizon variant (Appendix~\ref{app:long}, Proposition~\ref{prop:long}), and the control-theoretic verification (Appendix~\ref{app:control}).

\section{Main results}\label{sec:main}

Let $s>3/2$, $r=s-1>1/2$, and let $H^r_{\mathrm{ev}}(\T)$ denote the real, even Sobolev space. Define the affine space of unit-mass laws and its tangent space by
\[
 X=\Big\{m\in H^r_{\rm ev}(\T):\int_\T m=1\Big\},\qquad
 X_0=\Big\{h\in H^r_{\rm ev}(\T):\int_\T h=0\Big\}.
\]
For $j\ge1$, put $q_j=2\pi j$, $e_j(z)=\cos(q_jz)$ and $c_j(h)=2\int_\T h(z)e_j(z)\dd z$. On $X_0$ use the equivalent Hilbert norm
\[
 \norm{h}_*^2=\sum_{j\ge1}(1+q_j^2)^r c_j(h)^2=2\norm{h}_{H^r}^2.
\]
Write $P_d h=\sum_{j=1}^d c_j(h)e_j$, $Q_d=I-P_d$ and $\Es_d=Q_d X_0$. We abbreviate $P=P_1$, $Q=Q_1$ and $\Es=\Es_1$. The scaled product norm on $\R^d\times\Es_d$ is $\max\{|a|_\infty,\norm{y}_*\}$. Sobolev embedding gives $\norm{y}_{\infty}\le C_{s,d}\norm{y}_*$.

For a real, even running potential $V\in H^{s+2}(\T)$ and a strictly positive $m\in X$, solve on $[0,\tau]\times\T$
\begin{equation}\label{eq:PDE}
\begin{aligned}
 -u_t-\nu u_{zz}+\tfrac12|u_z|^2&=V(z),&u(\tau,\cdot)&=0,\\
 \rho_t&=\nu\rho_{zz}+\partial_z(\rho u_z),&\rho(0,\cdot)&=m.
\end{aligned}
\end{equation}
Set $A_V=\nu\partial_{zz}-V/(2\nu)$ and $S_V(t)=e^{tA_V}$. The terminal density and the frozen-law nonlinear update are
\begin{equation}\label{eq:Phi}
 \Phi_\tau(V,m)=S_V(\tau)\left(\frac{m}{S_V(\tau)1}\right),\qquad
 \Psi_{\varepsilon,\tau}(m)=\Phi_\tau(V_{\varepsilon,\tau}[m],m).
\end{equation}
By Proposition~\ref{prop:analytic} in Appendix~\ref{app:analytic}, $\Phi_\tau$ is positive and mass-preserving, locally analytic near $(0,1)$, and has relatively compact images on bounded small-potential neighbourhoods. The first variation is exactly
\begin{equation}\label{eq:variation}
D_m\Phi_\tau(0,1)h=S(\tau)h,\qquad
D_V\Phi_\tau(0,1)[e_j]=-\chi_j e_j,
\qquad S(t)=e^{\nu t\partial_{zz}},
\end{equation}
where $\chi_j$ is defined by \eqref{eq:introchi}; see Lemma~\ref{lem:susceptibility}.

Let $J\Subset U\subset\R$ with $J$ a nonempty compact interval and $U$ open, and let $f\in C^k(U)$, $k\ge1$. For $a\in U$ define
\begin{equation}\label{eq:b}
 b_\tau(a)=\frac{E_1a-f(a)}{\chi_1},\qquad
 V_{\varepsilon,\tau}[m]=\varepsilon b_\tau\!\left(\frac{c_1(m)}{\varepsilon}\right)e_1.
\end{equation}
This is the feedback \eqref{eq:introfb} written with $b_\tau$, and it is defined whenever $c_1(m)/\varepsilon\in U$. It is also a potential interaction with respect to the population: on each connected component of $U$ choose $a_0$ in it and let
\begin{equation}\label{eq:energy}
 \mathcal E_{\varepsilon,\tau}(m)
 =\frac12\int_{\varepsilon a_0}^{c_1(m)}
       \varepsilon b_\tau(z/\varepsilon)\dd z;
\qquad D\mathcal E_{\varepsilon,\tau}(m)[h]=\int_\T V_{\varepsilon,\tau}[m]h.
\end{equation}
No assertion of a globally defined energy on disconnected domains is made.

On $\mathcal D=J\times\overline B_1(\Es)$ write
\begin{equation}\label{eq:coords}
 m_\varepsilon(a,y)=1+\varepsilon(ae_1+y),\qquad
 \Psi_{\varepsilon,\tau}(m_\varepsilon(a,y))
 =1+\varepsilon\big(F_{\varepsilon,\tau}(a,y)e_1+G_{\varepsilon,\tau}(a,y)\big),
\end{equation}
so $T_{\varepsilon,\tau}=(F_{\varepsilon,\tau},G_{\varepsilon,\tau})$ is a map to $\R\times\Es$. The symbol $C^k(\mathcal D)$ denotes suprema of derivatives up to order $k$ on this closed product, computed from an open neighbourhood; compactness of the stable ball is neither assumed nor used.

\begin{theorem}[Fixed-horizon scalar realization in $C^k$]\label{thm:scalar}
Fix $\nu,\tau>0$, $s>3/2$, $J\Subset U$ and $f\in C^k(U)$ with $k\ge1$ finite. There are $C,\varepsilon_0>0$, depending only on $\nu,\tau,s,k,J,U$ and on the $C^k$ norm of $f$ on a compact neighbourhood of $J$ in $U$, such that for $0<\varepsilon<\varepsilon_0$ all $m_\varepsilon(a,y)$ are strictly positive probability densities, $T_{\varepsilon,\tau}$ is $C^k$ on a neighbourhood of $\mathcal D$, the image of the tube under $\Psi_{\varepsilon,\tau}$ is relatively compact in $H^r$, and
\begin{equation}\label{eq:mainestimate}
 \norm{T_{\varepsilon,\tau}-(f,S(\tau)|_{\Es})}_{C^k(\mathcal D)}\le C\varepsilon.
\end{equation}
The estimate is uniform over the entire unit stable ball and holds for every separately prescribed $\tau>0$.
\end{theorem}
\begin{proof}
For $(a,y)\in\mathcal D$ put $h=\varepsilon(ae_1+y)$ and $V=\varepsilon b_\tau(a)e_1$. Sobolev embedding gives $m_\varepsilon\ge1-C\varepsilon>0$. Because $J\Subset U$, the maps $(a,y)\mapsto(V,h)$ and each of their derivatives of order at most $k$ are $O(\varepsilon)$ uniformly on a slightly larger product neighbourhood. The analytic Taylor remainder from Proposition~\ref{prop:analytic} has the form
\[
\Phi_\tau(V,1+h)=1+S(\tau)h+D_V\Phi_\tau(0,1)[V]+R(V,h),
\]
where $R(0,0)=DR(0,0)=0$ and, on a fixed small open Banach-space ball, $\norm{R}=O(\norm{(V,h)}^2)$. Applying the higher-order chain rule to $R\circ(V,h)$ shows that all derivatives through order $k$ are $O(\varepsilon^2)$: in each term either a factor $DR=O(\varepsilon)$ multiplies an $O(\varepsilon)$ derivative, or at least two derivatives of $(V,h)$ occur (see Remark~\ref{rem:taylor}). This uses a bounded derivative of $\Phi_\tau$ of order at least $\max\{2,k\}$ and is valid on the full stable ball. The linear terms, by \eqref{eq:variation}, equal
\[
1+\varepsilon\{[E_1a-\chi_1b_\tau(a)]e_1+S(\tau)y\}
=1+\varepsilon\{f(a)e_1+S(\tau)y\}.
\]
Projecting onto the first mode and its complement and dividing by $\varepsilon$ gives \eqref{eq:mainestimate}. The positivity and compactness statements follow from the same local uniform analytic and smoothing estimates.
\end{proof}

The next theorem turns the scalar realization estimate into a statement about actual invariant sets of the parabolic density map, not merely of the approximate scalar map. Under strict expanding-covering hypotheses on $f$, it gives a compact invariant set of the exact response that factors onto the full one-sided $N$-shift, with periodic densities of every least period. The infinite-dimensional covering lemma, including the compact-core and periodic-closing steps, is proved in Appendix~\ref{app:cover}. The proof has three steps. (i) By \eqref{eq:mainestimate} with $k=1$, the exact map $T_{\varepsilon,\tau}$ lies $C^1$-close to $T_0(a,y)=(f(a),S(\tau)y)$ on the whole tube $\mathcal N=\bigcup_i I_i\times\overline B_1(\Es)$. (ii) Hence expansion, strict endpoint covering, stable contraction and the graph-cone inequalities hold for $T_{\varepsilon,\tau}$ itself, with explicit constants. (iii) Together with relative compactness of the image, Lemma~\ref{lem:cover} produces a compact invariant set of the exact response and its shift factor.

\begin{theorem}[Full symbolic shifts for the exact response]\label{thm:symbolic}
Suppose $f\in C^1(U)$ admits pairwise disjoint compact intervals $I_1,\ldots,I_N\subset J$, $N\ge2$, on each of which $f$ is monotone and
\[
\inf_{\cup_i I_i}|f'|>1.
\]
Assume that every endpoint image $f(I_i)$ covers the convex hull $[a_-,a_+]=\operatorname{conv}(\cup_iI_i)$ with a common strictly positive margin. For each fixed $\nu,\tau>0$ and all sufficiently small $\varepsilon>0$, the map $\Psi_{\varepsilon,\tau}$ of \eqref{eq:Phi}--\eqref{eq:b} has a nonempty compact set $\Lambda_\varepsilon$ of strictly positive densities such that
\[
\Psi_{\varepsilon,\tau}(\Lambda_\varepsilon)=\Lambda_\varepsilon,
\qquad \pi\circ\Psi_{\varepsilon,\tau}=\sigma\circ\pi
\]
for a continuous surjection $\pi$ onto the one-sided full $N$-shift, with shift map $\sigma$. Consequently $h_{\rm top}(\Psi_{\varepsilon,\tau}|_{\Lambda_\varepsilon})\ge\log N$ on this compact restriction, and periodic densities of every least period exist.
\end{theorem}
\begin{proof}
In scaled coordinates the limit is $T_0(a,y)=(f(a),S(\tau)y)$; its stable operator has norm
\[
 \alpha_0=\norm{S(\tau)|_{\Es}}=e^{-\nu q_2^2\tau}=e^{-\nu(4\pi)^2\tau}<1.
\]
On $\mathcal N=\bigcup_i(I_i\times\overline B_1(\Es))$, Theorem~\ref{thm:scalar} bounds the $C^1$ difference between $T_\varepsilon$ and $T_0$ by $\eta=C\varepsilon$. In the notation of Lemma~\ref{lem:cover}, $\lambda_u\ge\inf |f'|-\eta$, $\beta,\gamma\le\eta$ and $\alpha\le\alpha_0+\eta$. The slope $L=1$ satisfies the graph-cone inequalities for small $\eta$, since $\inf|f'|>1>\alpha_0$; stable trapping and strict endpoint covering also persist. Quantitatively, write $\lambda_f=\inf_{\cup I_i}|f'|$ and let $\mu$ be the smallest strict endpoint margin. First choose
\[
0<\eta<\min\left\{\frac\mu2,\frac{\lambda_f-1}{8},
                       \frac{1-\alpha_0}{8},\frac1{16}\right\},
\]
and then take $\varepsilon$ so that $C\varepsilon\le\eta$ and the entire Sobolev tube consists of positive densities. With unit graph slope, $\lambda\ge\lambda_f-2\eta>1$ and the stable graph contraction factor is bounded by $\alpha_0+\eta+(\alpha_0+2\eta)\eta/(\lambda_f-2\eta)<1$. The image of $\mathcal N$ is relatively compact by Theorem~\ref{thm:scalar}. Lemma~\ref{lem:cover} produces the compact shift-factor core and the periodic points. The coordinate map $(a,y)\mapsto1+\varepsilon(ae_1+y)$ transports the conclusion into $X$.
\end{proof}

\begin{corollary}[A two-shift cubic and unbounded finite entropies]\label{cor:examples}
For $f(a)=5a-5a^3$ choose $I_1=[-1/4,1/4]$ and $I_2=[3/4,21/20]$. The potential interaction
\begin{equation}\label{eq:cubic}
 V_{\varepsilon,\tau}[m]
 =\left(-\frac{5-E_1}{\chi_1}c_1(m)
       +\frac{5}{\chi_1\varepsilon^2}c_1(m)^3\right)e_1
\end{equation}
produces a full two-shift factor for every fixed $\nu,\tau>0$ and small $\varepsilon$. More generally, for each integer $N\ge2$ there is a polynomial rank-one potential family whose exact response has entropy at least $\log N$ on a compact invariant set contained in any preassigned $H^r$ neighbourhood of the uniform density $1$. The feedback family depends on $N$ and on that neighbourhood.
\end{corollary}
\begin{proof}
For the cubic, the design \eqref{eq:b} reduces to \eqref{eq:cubic}, and the branch endpoint values are $-75/64$, $75/64$, $105/64$ and $-861/1600$, respectively. They strictly cover $[-1/4,21/20]$; moreover, the minimum absolute derivatives of $f$ on the two intervals are $65/16$ and $55/16$. Thus Theorem~\ref{thm:symbolic} applies.
For general $N$, let $T_N(\cos\theta)=\cos(N\theta)$ be the degree-$N$ Chebyshev polynomial. In each of its $N$ monotonicity laps, choose an interior compact interval on which $|T_N'|$ is bounded below by a positive constant and whose endpoint values straddle $[-a,a]$ for a fixed $a\in(0,1)$. These $N$ intervals are disjoint and lie in $[-1,1]$. For sufficiently large fixed $c>0$, the profile $f=cT_N$ has derivative greater than one in absolute value on all intervals, and each branch image strictly covers $[-1,1]$. Apply Theorem~\ref{thm:symbolic}; finally decrease $\varepsilon$ to place the whole scaled density tube inside the prescribed neighbourhood of $1$.
\end{proof}

The next result extends the scalar construction to several observable modes. It does not rely on invariance of a finite-dimensional Fourier manifold; no such exact invariant manifold is assumed or needed.

\begin{theorem}[Full-space $C^k$ realization of finite-mode maps at a prescribed horizon]\label{thm:multimode}
Fix $d\ge1$, $\nu,\tau>0$, $s>3/2$, a compact set $K\Subset U\subset\R^d$, and a map $F=(F_1,\ldots,F_d)\in C^k(U;\R^d)$ for finite $k\ge1$. For $a_j=c_j(m)/\varepsilon$ define
\begin{equation}\label{eq:multifb}
 V^{(d)}_{\varepsilon,\tau}[m]
 =\varepsilon\sum_{j=1}^d\frac{E_ja_j-F_j(a)}{\chi_j}\,e_j.
\end{equation}
On $K\times\overline B_1(\Es_d)$ write $T^{(d)}_{\varepsilon,\tau}=(\mathfrak F_\varepsilon,\mathfrak G_\varepsilon)$ for the exact response in scaled coordinates, for $d=1$ and $F=f$ this is the map $T_{\varepsilon,\tau}$ of \eqref{eq:coords}; explicitly,
\[
 \Phi_\tau\big(V^{(d)}_{\varepsilon,\tau}[m],m\big)
 =1+\varepsilon\left(\sum_{j=1}^d\mathfrak F_{\varepsilon,j}(a,y)e_j+\mathfrak G_\varepsilon(a,y)\right),
 \quad m=1+\varepsilon\left(\sum_{j=1}^da_je_j+y\right).
\]
For every $\tau>0$ there are $C=C_\tau$ and $\varepsilon_0(\tau)>0$ such that for $0<\varepsilon<\varepsilon_0(\tau)$ these densities are positive, the image of the tube is relatively compact, and
\begin{equation}\label{eq:multimain}
\norm{T^{(d)}_{\varepsilon,\tau}
 -(F,S(\tau)|_{\Es_d})}_{C^k(K\times\overline B_1(\Es_d))}
 \le C\varepsilon,
\qquad \norm{S(\tau)|_{\Es_d}}=e^{-\nu q_{d+1}^2\tau}<1.
\end{equation}
For $d>1$, the constructed $V^{(d)}_{\varepsilon,\tau}$ need not be the variational derivative of a scalar energy of the population (see Theorem~\ref{thm:potential}).
\end{theorem}
\begin{proof}
The map $(a,y)\mapsto(V^{(d)}[m],m-1)$, together with its derivatives through order $k$, is uniformly $O(\varepsilon)$ on the tube. Apply the Taylor-remainder and chain-rule argument of Theorem~\ref{thm:scalar}, now retaining modes $1,\ldots,d$ separately, which gives \eqref{eq:multimain}. The first-order coefficient in mode $j$ is $E_ja_j-\chi_j(E_ja_j-F_j(a))/\chi_j=F_j(a)$. Fourier orthogonality makes the other first-order terms equal to $S(\tau)y$. Sobolev embedding and positive-time smoothing prove the auxiliary claims.
\end{proof}

The realization estimate \eqref{eq:multimain} is an approximation statement for the full response. In the scalar expanding-covering setting of Theorem~\ref{thm:symbolic}, it yields an exact invariant-set conclusion for that response; the selected moments need not form an autonomous subsystem for this to hold (compare Proposition~\ref{prop:closure}).

The next theorem improves the existence of periodic points: for an individual nondegenerate target cycle it gives local uniqueness, an approximation rate, and information on the whole linearized infinite-dimensional response.

\begin{theorem}[Continuation and hyperbolicity of a prescribed cycle]\label{thm:cycles}
In Theorem~\ref{thm:scalar}, assume $f\in C^2(U)$ and let $a_0,\ldots,a_{n-1}\in\operatorname{int}J$ be a period-$n$ orbit of $f$, with least period $n$, such that
\[
 \lambda_*=(f^n)'(a_0)\ne1.
\]
Then, for small positive $\varepsilon$, there is a unique period-$n$ point $z_\varepsilon=(a_\varepsilon,y_\varepsilon)$ of the scaled response in a fixed sufficiently small neighbourhood of $(a_0,0)$, and
\[
 |a_\varepsilon-a_0|+\norm{y_\varepsilon}_*\le C_n\varepsilon.
\]
The corresponding density is $1+\varepsilon(a_\varepsilon e_1+y_\varepsilon)$; its response orbit has least period $n$. If $|\lambda_*|>1$, the linearization of $T_\varepsilon^n$ at $z_\varepsilon$ has exactly one unstable eigenvalue $\lambda_\varepsilon=\lambda_*+O(\varepsilon)$, of algebraic multiplicity one. If $|\lambda_*|<1$, the continued period-$n$ cycle is attracting. For $|\lambda_*|>1$, all remaining spectral values lie in
\[
 \big\{z\in\mathbb C:|z|\le e^{-\nu q_2^2n\tau}+C_n\varepsilon\big\}\subset\{ |z|<1\}.
\]
\end{theorem}
\begin{proof}
Let $T_0(a,y)=(f(a),S(\tau)y)$. Since the chosen scalar cycle is finite and lies in the interior of the coordinate domain, select small neighbourhoods $U_j$ of its successive points such that all first $n$ iterates of $T_0$ are defined on a neighbourhood $W$ of $(a_0,0)$. Uniform $C^1$ convergence on these neighbourhoods gives the same property for $T_\varepsilon$ when $\varepsilon$ is small. No forward invariance of the entire interval $J$ is needed.

Put $g(a,y)=(b_\tau(a)e_1,ae_1+y)$ and $\widehat\Phi_\tau(V,h)=\Phi_\tau(V,1+h)$. The scaled difference quotient has the removable singularity
\[
 T_\varepsilon(a,y)=\mathcal C\int_0^1
 D\widehat\Phi_\tau(t\varepsilon g(a,y))[g(a,y)]\,\dd t,
\]
where $\mathcal C$ is the coefficient/stable projection. This extends the scaled response to a jointly $C^1$ family in $(\varepsilon,a,y)$, with $C^2$ dependence on the state under the stated regularity. Writing $H(\varepsilon,z)=T_\varepsilon^n(z)-z$ gives
\[
 D_zH(0,(a_0,0))=\operatorname{diag}
 (\lambda_*-1,S(n\tau)|_{\Es}-I).
\]
The diagonal blocks are invertible, and the norm of the inverse in the product coordinates is at most
\[
 M_n=\max\left\{\frac1{|\lambda_*-1|},
       \frac1{1-e^{-\nu q_2^2n\tau}}\right\}.
\]
The implicit-function theorem supplies the unique local continuation and, since $H(\varepsilon,(a_0,0))=O(\varepsilon)$, a displacement bounded by $C_n\varepsilon$. For instance, on a small ball where the normalized Newton derivative differs from the identity by at most $1/2$, the contraction bound is $2M_n\|H(\varepsilon,(a_0,0))\|$. Separation of the original cycle points preserves the least period.

For the spectral claims, use the equivalent complex Hilbert norm
\[
 \|(a,y)\|_{\mathcal H}^2=|a|^2+\|y\|_*^2,
 \qquad \mathcal H=\mathbb C\oplus\Es_{\mathbb C}.
\]
In this norm the limiting derivative $L_0=\operatorname{diag}(\lambda_*,S(n\tau)|_{\Es})$ is normal, and
\[
 L_\varepsilon=DT_\varepsilon^n(z_\varepsilon)
                  =L_0+E_\varepsilon,
 \qquad\|E_\varepsilon\|_{\mathcal L(\mathcal H)}\le C_n\varepsilon.
\]
Normality is used only for $L_0$, for which $\|(z-L_0)^{-1}\|=\operatorname{dist}(z,\sigma(L_0))^{-1}$. The resolvent Neumann series confines $\sigma(L_\varepsilon)$ to the $C_n\varepsilon$-neighbourhood of $\sigma(L_0)$. If $|\lambda_*|>1$, a contour separated from the stable spectrum encloses a rank-one Riesz cluster with a real eigenvalue $\lambda_\varepsilon=\lambda_*+O(\varepsilon)$. The remaining spectrum is within the asserted stable disk, proving the expansion statement. If instead $|\lambda_*|<1$, all multipliers have modulus strictly below one, so the continued cycle is attracting. The excluded case $\lambda_*=1$ requires other methods; the case $\lambda_*=-1$ is permitted for continuation but does not imply hyperbolicity.
\end{proof}

\begin{remark}[Nondegeneracy does not imply hyperbolicity]\label{rem:nonhyp}
Take the target $f(a)=-a$ and its fixed point $a_0=0$. Then $\lambda_*=f'(0)=-1\ne1$, but the designed feedback simplifies to
\[
 V_{\varepsilon,\tau}[m]=\frac{E_1+1}{\chi_1}c_1(m)e_1,
\]
independently of $\varepsilon$. The uniform density is an exact fixed point, and its first Fourier multiplier is
\[
 D\Psi_{\varepsilon,\tau}(1)e_1
 =\left(E_1-\chi_1\frac{E_1+1}{\chi_1}\right)e_1=-e_1.
\]
Thus the fixed point is nonhyperbolic at every amplitude, although it satisfies the nondegeneracy condition used for continuation. The expanding condition $|\lambda_*|>1$ in Theorem~\ref{thm:cycles} is essential for its hyperbolicity conclusion.
\end{remark}

\begin{theorem}[Robustness to scaled interaction errors]\label{thm:robust}
Assume the expanding covering hypotheses of Theorem~\ref{thm:symbolic}, and put $\mathcal N=\bigcup_i(I_i\times\overline B_1(\Es))$. There exist $\delta_0>0$ and $\varepsilon_0>0$ with the following property. Let $\widetilde V_\varepsilon[m]$ be any real even potential, $C^1$ in the scaled tube variables on an open neighbourhood of $\mathcal N$, such that
\begin{equation}\label{eq:perturb}
 \norm{\widetilde V_\varepsilon[m_\varepsilon(\cdot)]
  -V_{\varepsilon,\tau}[m_\varepsilon(\cdot)]}
 _{C^1(\mathcal N;H^{s+2})}\le\varepsilon\delta_0.
\end{equation}
For every $0<\varepsilon<\varepsilon_0$, the exact PDE response with $\widetilde V_\varepsilon$ has a compact strictly invariant set of positive densities factoring onto the full one-sided $N$-shift, and has periodic points of all least periods. No variational or convolution structure is imposed on the perturbation.
\end{theorem}
\begin{proof}
Write $\Delta V=\widetilde V_\varepsilon-V_{\varepsilon,\tau}$, so that \eqref{eq:perturb} bounds $\Delta V$ by $\varepsilon\delta_0$ in $C^1$. Local analyticity of $\Phi_\tau$ and uniform boundedness of its first two derivatives on the small potential/density ball imply, by the mean-value formula and the chain rule,
\[
\norm{\Phi_\tau(\widetilde V_\varepsilon[m_\varepsilon(\cdot)],m_\varepsilon(\cdot))
-\Phi_\tau(V_{\varepsilon,\tau}[m_\varepsilon(\cdot)],m_\varepsilon(\cdot))}
_{C^1(\mathcal N;H^r)}\le C\varepsilon\delta_0.
\]
After division by $\varepsilon$, the scaled response differs by at most $C\delta_0$ in $C^1$, so its distance to $T_0$ is at most $C(\varepsilon+\delta_0)$. Choose $\delta_0$ and then $\varepsilon_0$ so that the strict derivative inequalities, strict endpoint crossings and stable trapping of the covering lemma are retained. Positive-time parabolic smoothing and the uniform smallness of $\widetilde V_\varepsilon$ ensure relatively compact images of the relevant stable tube. Apply Lemma~\ref{lem:cover}.
\end{proof}

\begin{theorem}[A general cost of nonaffine microscopic realization]\label{thm:cost}
Fix $s>3/2$, $\nu,\tau>0$ and a compact interval $J$. Let $V_\varepsilon$ be a family of arbitrary $C^2$ interactions defined on open neighbourhoods of the scaled segment
$\{1+\varepsilon ae_1:a\in J\}$, with values in a fixed sufficiently small ball of $H^{s+2}_{\rm ev}$. Put
\[
 M_{1,\varepsilon}=\sup_{a\in J}
 \|D_mV_\varepsilon[1+\varepsilon ae_1]\|_{\mathcal L(X_0,H^{s+2}_{\rm ev})},
 \qquad M_{2,\varepsilon}=\sup_{a\in J}
 \|D_m^2V_\varepsilon[1+\varepsilon ae_1]\|_{\mathcal L^2(X_0;H^{s+2}_{\rm ev})},
\]
where $\mathcal L^2(X_0;H^{s+2}_{\rm ev})$ denotes bounded bilinear maps on $X_0\times X_0$ and all norms on $X_0$ are inherited from $H^r$. Only segment suprema are needed. Suppose that
\[
 \mathfrak f_\varepsilon(a)=\frac1\varepsilon c_1\bigl(\Phi_\tau(V_\varepsilon[1+\varepsilon ae_1],1+\varepsilon ae_1)\bigr)
\]
converges in $C^1(J)$ to a nonaffine $f\in C^1(J)$. Then, for constants $c,C>0$ independent of sufficiently small $\varepsilon$,
\begin{equation}\label{eq:universal-cost}
 M_{2,\varepsilon}+(1+M_{1,\varepsilon})^2\ge \frac c\varepsilon.
\end{equation}
If $M_{1,\varepsilon}$ is uniformly bounded, then $M_{2,\varepsilon}\ge c\varepsilon^{-1}-C$. Conversely, let $f\in C^2(U)$ with $f''\not\equiv0$ on $J$. Then the rank-one design \eqref{eq:b} of Theorem~\ref{thm:scalar} has uniformly bounded $M_{1,\varepsilon}$ and satisfies $M_{2,\varepsilon}\asymp\varepsilon^{-1}$ on the segment; hence the order is attained in that class.
\end{theorem}
\begin{proof}
The analytic response has uniformly bounded first and second Fr\'echet derivatives on the small potential--density ball specified above (see Appendix~\ref{app:analytic}). Write $\Psi_\varepsilon(m)=\Phi_\tau(V_\varepsilon[m],m)$. The exact second chain rule is
\begin{align*}
D^2\Psi_\varepsilon(m)[h,k]
 &=D^2\widehat\Phi_\tau(V_\varepsilon[m],m-1)
       [(DV_\varepsilon[m]h,h),(DV_\varepsilon[m]k,k)]\\
 &\quad+D_V\Phi_\tau(V_\varepsilon[m],m)[D^2V_\varepsilon[m](h,k)].
\end{align*}
Consequently $\|D^2\Psi_\varepsilon\|\le C[(1+M_{1,\varepsilon})^2+M_{2,\varepsilon}]$ along the segment. Since $f$ is not affine there are $a,b\in J$ with $\Delta=|f'(a)-f'(b)|>0$. Set $\eta_\varepsilon=\|\mathfrak f_\varepsilon-f\|_{C^1(J)}$. Differentiation along the segment gives
\[
 \mathfrak f_\varepsilon''(t)
 =\varepsilon c_1\left(D^2\Psi_\varepsilon(1+\varepsilon te_1)[e_1,e_1]\right),
\]
whence
\[
\Delta-2\eta_\varepsilon\le C\varepsilon |a-b|
 \big((1+M_{1,\varepsilon})^2+M_{2,\varepsilon}\big).
\]
For small $\varepsilon$ we have $\eta_\varepsilon<\Delta/4$, which gives \eqref{eq:universal-cost}. For the designed scalar feedback,
\[
DV_{\varepsilon,\tau}[m]h=b_\tau'(c_1(m)/\varepsilon)c_1(h)e_1,
\qquad
D^2V_{\varepsilon,\tau}[m](h,k)
 =\frac{b_\tau''(c_1(m)/\varepsilon)}\varepsilon c_1(h)c_1(k)e_1.
\]
These identities give the upper bounds $M_{1,\varepsilon}=O(1)$ and $M_{2,\varepsilon}=O(\varepsilon^{-1})$. For the lower bound, note that $b_\tau''=-f''/\chi_1$, so $\sup_J|b_\tau''|>0$ because $f''\not\equiv0$ on $J$. In the star norm the exact bilinear operator norm at the point $a$ is $|f''(a)|\,\|e_1\|_{H^{s+2}}/(\chi_1\varepsilon\|e_1\|_*^2)$, and taking the supremum over $a\in J$ gives $M_{2,\varepsilon}\ge c\varepsilon^{-1}$.
\end{proof}

\begin{remark}[Quantifiers and interpretation]\label{rem:cost}
Here $X_0$ is the tangent space of the mass-one affine space $X$, with the inherited $H^r$ norm. Derivatives of an interaction are required only in directions of $X_0$; no extension to an open subset of the whole $H^r$ space is used. The bound \eqref{eq:universal-cost} does not assert $M_{2,\varepsilon}\gtrsim\varepsilon^{-1}$ when $M_{1,\varepsilon}$ is unbounded: the first derivative may itself be costly. It does not prohibit complicated finite-amplitude behaviour of a fixed convolution model. It does show why a family uniformly bounded in both population derivatives cannot realize a genuinely nonlinear scalar map on shrinking $\varepsilon$-tubes. No assumption from this theorem enters the positive realization results.
\end{remark}

\begin{theorem}[Potentiality and the leading finite-mode targets]\label{thm:potential}
In Theorem~\ref{thm:multimode}, let $U\subset\R^d$ be an open coordinate domain and $F\in C^1(U;\R^d)$. Set $D_\chi=\operatorname{diag}(\chi_1,\ldots,\chi_d)$ and $E=\operatorname{diag}(E_1,\ldots,E_d)$. Define the affine moment-saturated domain
\[
 \Omega_\varepsilon=\bigl\{m\in X:\ c(m)/\varepsilon\in U\bigr\},
 \qquad c(m)=(c_1(m),\ldots,c_d(m)).
\]
The designed interaction \eqref{eq:multifb} is the first population variation of a $C^1$ energy on $\Omega_\varepsilon$ (with derivatives evaluated on $X_0$) if and only if the one-form
\[
 \omega=b(a)\cdot\dd a,\qquad b(a)=D_\chi^{-1}(Ea-F(a)),
\]
is exact on $U$. Locally this is equivalent to
\begin{equation}\label{eq:weighted}
D_\chi^{-1}DF(a)\ \text{being symmetric},
\quad\text{i.e.}\quad
\frac{\partial_iF_j(a)}{\chi_j}=\frac{\partial_jF_i(a)}{\chi_i}.
\end{equation}
On any convex coordinate domain, \eqref{eq:weighted} is also sufficient globally. On a general domain, global exactness, and not closedness alone, is required. Wherever the weighted symmetry condition holds, $DF(a)$ has only real eigenvalues. If the energy is defined solely for strictly positive densities, the same criterion applies on a coordinate subdomain admitting a positive affine section (in particular, any fixed bounded neighbourhood of the compact target set for all sufficiently small $\varepsilon$).
\end{theorem}
\begin{proof}
For $m\in\Omega_\varepsilon$ and $h\in X_0$, the prescribed feedback is
\[
 V_\varepsilon[m]=\varepsilon\sum_{j=1}^d b_j(c(m)/\varepsilon)e_j,
 \qquad D\mathcal E_\varepsilon(m)[h]=\int_\T V_\varepsilon[m]h.
\]
Necessity of exactness does not require the energy to depend only on the moments. Indeed, the affine section
\[
 \sigma_\varepsilon(a)=1+\varepsilon\sum_{j=1}^d a_je_j,
 \qquad a\in U,
\]
satisfies $c(\sigma_\varepsilon(a))=\varepsilon a$ and hence lies in $\Omega_\varepsilon$. If such a $C^1$ energy exists, its restriction gives
\[
 \Gamma(a)=\frac{2}{\varepsilon^2}
                \mathcal E_\varepsilon(\sigma_\varepsilon(a)).
\]
Since $\partial_j\sigma_\varepsilon=\varepsilon e_j$ and $\int_\T e_je_k=\delta_{jk}/2$, the chain rule yields
\[
 \partial_j\Gamma(a)
 =\frac{2}{\varepsilon^2}
     D\mathcal E_\varepsilon(\sigma_\varepsilon(a))[\varepsilon e_j]
 =2\sum_{k=1}^d b_k(a)\int_\T e_ke_j=b_j(a).
\]
For $d=1$ the energy \eqref{eq:energy} is of this form. Thus $\omega=d\Gamma$, whether or not the original energy had a moment-only representation.

Conversely, if $\omega$ is exact on $U$, choose a primitive $\Gamma\in C^2(U)$ with $\nabla\Gamma=b$ and define
\[
 \mathcal E_\varepsilon(m)=\frac{\varepsilon^2}{2}
             \Gamma(c(m)/\varepsilon),\qquad m\in\Omega_\varepsilon.
\]
For every $h\in X_0$,
\[
 D\mathcal E_\varepsilon(m)[h]
 =\varepsilon\sum_{j=1}^d (\partial_j\Gamma)(a)
                         \int_\T he_j
 =\int_\T V_\varepsilon[m]h.
\]
This proves sufficiency on the whole affine domain, with no positivity or convexity assumption. The differential criterion is the familiar exact-one-form test \cite{MS,Letcher}; here, differentiation of $b=D_\chi^{-1}(Ea-F)$ gives
\[
 Db=D_\chi^{-1}E-D_\chi^{-1}DF.
\]
Because the first term is diagonal, local closedness is equivalent to \eqref{eq:weighted}. On a convex $U$, a primitive can be written, for any $a_0\in U$, as
\[
 \Gamma(a)=\int_0^1 b(a_0+t(a-a_0))\cdot(a-a_0)\,\dd t.
\]
The criterion is only local when $U$ is not simply connected: a closed form need not be exact. When $\omega$ is locally exact,
\[
 D_\chi^{-1/2}DF(a)D_\chi^{1/2}
 =E-D_\chi^{1/2}D^2\Gamma(a)D_\chi^{1/2},
\]
which is symmetric and similar to $DF(a)$, giving the real-spectrum conclusion. Finally, if energies are defined only for positive laws, restrict the necessity argument to a coordinate subdomain where $\sigma_\varepsilon(a)$ is positive. Each bounded subdomain has this property for sufficiently small $\varepsilon$ by Sobolev embedding; the sufficiency construction restricts to the positive part of $\Omega_\varepsilon$.
\end{proof}

\begin{remark}[A concrete obstruction]
For $d=2$, the target $F(a_1,a_2)=(-a_2,a_1)$ fails \eqref{eq:weighted}, since it would require $-1/\chi_1=1/\chi_2$. It is therefore achievable by the nonpotential finite-mode inverse formula but not by a scalar population energy of the specified moment form. This is not a classification of all finite-amplitude law responses or of all population interactions.
\end{remark}

\begin{proposition}[Failure of exact scalar closure in the rank-one family]\label{prop:closure}
For a fixed horizon $\tau>0$, let $\mathcal L_vh=S_{ve_1}(\tau)[h/S_{ve_1}(\tau)1]$ and $a_1=\nu q_1^2$, $E_1=e^{-a_1\tau}$. Then
\begin{equation}\label{eq:closure-coef}
\left.\partial_v c_1(\mathcal L_ve_2)\right|_{v=0}
=\frac{E_1(1-E_1)^2(E_1+2)}{12\nu a_1}>0.
\end{equation}
Consequently, for the response with feedback \eqref{eq:b}, uniformly for $a\in J$ and $\|y\|_*\le1$,
\[
D_yF_{\varepsilon,\tau}(a,y)e_2
=\varepsilon b_\tau(a)
\frac{E_1(1-E_1)^2(E_1+2)}{12\nu a_1}+O(\varepsilon^2b_\tau(a)^2).
\]
In particular, if $b_\tau\not\equiv0$ on $J$, there are $a\in J$ and $\varepsilon_0>0$ such that $D_yF_{\varepsilon,\tau}(a,\cdot)e_2\ne0$ for all $0<\varepsilon<\varepsilon_0$, so the first Fourier moment is not an exactly autonomous scalar factor on the tube. In the exceptional case $b_\tau\equiv0$, equivalently $f(a)=E_1a$, the potential vanishes and exact scalar closure holds: $F_{\varepsilon,\tau}(a,y)=E_1a$.
\end{proposition}
\begin{proof}
The rank-one feedback is constant along the $e_2$ direction, so differentiating the exact response with respect to the stable input gives $D_yF_{\varepsilon,\tau}(a,y)e_2=c_1(\mathcal L_{\varepsilon b_\tau(a)}e_2)$. The $v$-derivative of the semigroup is given by Duhamel's formula. Since $e_1e_2=(e_1+e_3)/2$,
\[
 c_1\left(\left.\partial_vS_{ve_1}(\tau)e_2\right|_{v=0}\right)
 =-\frac{E_1(1-E_1^3)}{12\nu a_1}.
\]
On the other hand, $\left.\partial_v[S_{ve_1}(\tau)1]^{-1}\right|_{v=0}=(1-E_1)e_1/(2\nu a_1)$. Applying $S(\tau)$ to its product with $e_2$ contributes $E_1(1-E_1)/(4\nu a_1)$ to the first coefficient. Adding these terms yields \eqref{eq:closure-coef}. Analyticity in $v$ gives the stated expansion with remainder $O(v^2)$, and $v=\varepsilon b_\tau(a)$ with $b_\tau$ bounded on the compact set $J$. If $b_\tau(a)\ne0$ at some $a\in J$, the leading term is nonzero and dominates the remainder once $\varepsilon|b_\tau(a)|$ is small, so $D_yF_{\varepsilon,\tau}(a,\cdot)e_2\ne0$. Conversely, if $b_\tau\equiv0$, then $V_{\varepsilon,\tau}=0$ and the Fourier heat flow gives $F_{\varepsilon,\tau}(a,y)=E_1a$ exactly. This does not contradict the $O(\varepsilon)$ approximate realization.
\end{proof}

\section{Conclusions and open problems}\label{sec:conclusion}
We have used an exact finite-horizon semigroup formula to isolate an invertible Fourier susceptibility and to program prescribed leading observable dynamics into an infinite-dimensional positive law map. Strict covering produces compact invariant shift factors of the exact PDE response, so the conclusions do not concern a Galerkin approximation. Nondegenerate individual cycles continue, but hyperbolicity requires a separate strict multiplier hypothesis. The construction is local in amplitude and depends on the target, the amplitude, the number of modes and the horizon; none of the estimates is uniform as $\tau\downarrow0$, $d\to\infty$, or as the alphabet grows. Indeed, $\chi_j\sim q_j^2\tau^2/2$ as $\tau\downarrow0$, so even the linear inverse design becomes increasingly expensive at short horizons.

The derivative-cost theorem makes a standard second-chain-rule obstruction quantitative at the present microscopic scale and identifies the sharp order $\varepsilon^{-1}$ when the first population derivative is bounded. The potentiality result specializes the classical exact-one-form test \cite{MS,Letcher} to the explicitly computed susceptibilities; neither the integrability principle nor its mixed-partial criterion is claimed as new. The closure calculation supplies a model-specific nonzero coefficient, showing why robust one-step approximation with nontrivial feedback is not exact scalar semiconjugacy; the zero-feedback heat response is the explicit exception.

Feew questions remain.
\begin{enumerate}
\item Which finite-amplitude law maps can be realized subject to a uniform regularity budget, convex population energies, or translation-invariant convolutions, without writing the desired nonlinearity into the interaction?
\item Invariant smooth manifolds or exact observable closure generally require feedback depending on stable modes, which changes the admissible class. For nonhyperbolic targets, approximation alone cannot decide bifurcations or exact Cantor conjugacies.
\item The Fokker--Planck response has a controlled-diffusion interpretation, but a finite-player or empirical-measure limit in a Wasserstein-type topology requires estimates beyond the positive Sobolev-density setting used here.
\end{enumerate}

\appendix
\section{Positive semigroups, analytic dependence and susceptibility}\label{app:analytic}

We record the standard parabolic facts used above; references for analytic semigroups, smoothing and parameter-dependent spectral theory include \cite{Davies,Henry,Kato}. They are collected here to separate analytic input from the engineered feedback constructions.

\begin{proposition}[Exact formula, positivity and local compactness]\label{prop:analytic}
For fixed $\nu,\tau>0$ and every real even potential $V\in H^{s+2}(\T)$, the solution of \eqref{eq:PDE} has terminal law \eqref{eq:Phi}. For each fixed $V$, the map $m\mapsto\Phi_\tau(V,m)$ is linear, positive and mass preserving. The map $(V,m)\mapsto\Phi_\tau(V,m)$ is analytic on a sufficiently small neighbourhood of $(0,1)$ in $H^{s+2}_{\rm ev}\times X$ (with strictly positive $m$), with values in $X$. For fixed $\tau>0$, every uniformly bounded small neighbourhood in these variables has image bounded in $H^{r+\theta}$ for some $\theta>0$; its image is therefore relatively compact in $H^r$.
\end{proposition}
\begin{proof}
Set $z=\exp(-u/(2\nu))$ and $v=\rho/z$. A direct substitution into \eqref{eq:PDE} yields $z_t=-A_Vz$, $z(\tau)=1$, and $v_t=A_Vv$, $v(0)=m/z(0)$. Thus $z(t)=S_V(\tau-t)1$ and $v(t)=S_V(t)[m/S_V(\tau)1]$. At the terminal time $z(\tau)=1$, which proves the formula. The parabolic maximum principle or the Feynman--Kac formula implies $S_V(\tau)1>0$. Since $A_V$ is self-adjoint, for all unit-mass $m$
\[
 \int_\T \Phi_\tau(V,m)=\ip{S_V(\tau)1}{m/S_V(\tau)1}=\int_\T m=1.
\]
Strict positivity for positive $m$ follows from positivity of the heat kernel. We now give the uniform functional-analytic bounds behind these assertions. In $H^r$, multiplication by $V$ has operator norm at most $c_r\|V\|_{H^r}$ because $r>1/2$. Iterated Duhamel terms for $S_V(t)$ therefore satisfy
\[
\left\|\int_{0<t_1<\cdots<t_n<t}
 S(t-t_n)M_{-V/(2\nu)}\cdots
 M_{-V/(2\nu)}S(t_1)\,\dd\boldsymbol{t}\right\|_{H^r\to H^r}
\le\frac1{n!}\left(\frac{c_rt\|V\|_{H^r}}{2\nu}\right)^n.
\]
Consequently $V\mapsto S_V(t)$ is analytic in operator norm and $\|S_V(t)\|_{H^r\to H^r}\le e^{\kappa t}$, where $\kappa=c_r\|V\|_{H^r}/(2\nu)$. The same expansion in $H^{s+2}$ gives $Z_V=S_V(\tau)1=1+O(\|V\|_{H^{s+2}})$ in that Banach algebra. On a sufficiently small fixed ball, the convergent Neumann series for $(1+(Z_V-1))^{-1}$ yields uniform $H^{s+2}$ bounds and analytic dependence of $Z_V^{-1}$; a pointwise positive denominator alone would not give this. Thus the shifted map
\[
 \widehat\Phi_\tau(V,h)=S_V(\tau)[(1+h)Z_V^{-1}],
\qquad h\in X_0,
\]
is jointly analytic from a neighbourhood of $(0,0)$ into $H^r$ and takes values in the affine mass-one space.

For relative compactness, fix $0<\theta<2$. The heat estimate $\|S(t)g\|_{H^{r+\theta}}\le C_\theta[1+(\nu t)^{-\theta/2}]\|g\|_{H^r}$ and one Duhamel step give
\begin{align*}
\|S_V(\tau)g\|_{H^{r+\theta}}
&\le C_\theta[1+(\nu\tau)^{-\theta/2}]\|g\|_{H^r}\\
&\quad+C_\theta\kappa e^{\kappa\tau}
\left[\tau+\frac{\nu^{-\theta/2}\tau^{1-\theta/2}}{1-\theta/2}\right]\|g\|_{H^r}.
\end{align*}
The factor $VS_V(t)g$ is estimated in $H^r$, so no smoothing property of multiplication is invoked. Apply this with $g=(1+h)Z_V^{-1}$, uniformly bounded in $H^r$ on the indicated balls. The image is bounded in $H^{r+\theta}$, whose embedding into $H^r$ is compact. The closed infinite-dimensional input ball is not assumed compact. In terms of the PDE, $z(t)=S_V(\tau-t)1$ belongs to $C([0,\tau];H^{s+2})\cap C^1([0,\tau];H^s)$, and $u=-2\nu\log z$ has the same parabolic regularity. The Fokker--Planck solution is mild in $H^r$ at time zero and has finite positive-time spatial regularity: for each $0<\theta<2$,
\[
 \rho\in C([0,\tau];H^r_{\rm ev})
 \cap C((0,\tau];H^{r+\theta}_{\rm ev}),
\]
with $H^{r+\theta}$ bounds uniform on closed positive-time subintervals. To see this for intermediate times, use $\rho(t)=z(t)S_V(t)[m/Z_V]$ and the same heat/Duhamel estimate on $[\delta,\tau]$, noting that multiplication by $z(t)\in H^{s+2}$ preserves $H^{r+\theta}$ for $\theta<2$. Because $r>1/2$, we may choose $\theta<2$ with $r+\theta>5/2$, which gives enough spatial regularity for the forward equation whenever $t>0$. Positivity follows from the positive parabolic kernel. These statements do not imply $C^\infty$ smoothing for coefficients with only finite Sobolev regularity.
\end{proof}

\begin{lemma}[Diagonal susceptibility]\label{lem:susceptibility}
The derivatives of $\Phi_\tau$ at $(0,1)$ are given by \eqref{eq:variation}. For each mean-zero trigonometric polynomial $H$, and hence by continuity in the indicated Sobolev spaces,
\[
 D_V\Phi_\tau(0,1)[H]
 =-\frac1{2\nu^2}(-\partial_{zz})^{-1}(I-S(\tau))^2H.
\]
\end{lemma}
\begin{proof}
The derivative in $m$ follows from the exact formula with $V=0$. For $V=t e_j$, Duhamel's formula gives
\[
 \left.\frac\dd{\dd t}S_{t e_j}(\tau)1\right|_{t=0}
 =-\frac1{2\nu}\int_0^\tau S(\tau-z)e_j\dd z
 =-\frac{1-E_j}{2\nu^2q_j^2}e_j.
\]
Differentiation of the denominator and of the numerator in \eqref{eq:Phi} gives
\[
 D_V\Phi_\tau(0,1)[e_j]
 =-\frac{1-E_j}{2\nu^2q_j^2}(1-E_j)e_j=-\chi_j e_j.
\]
By linearity and Fourier expansion, the operator formula follows.
\end{proof}

\begin{remark}[Control of the Taylor remainder]\label{rem:taylor}
The second derivative of the analytic map $\Phi_\tau$ is uniformly bounded on a sufficiently small open Banach ball around $(0,1)$. Therefore, with $w=(V,h)$,
\[
\widehat\Phi_\tau(V,h)-1-S(\tau)h-D_V\Phi_\tau(0,1)[V]
 =\int_0^1(1-t)D^2\widehat\Phi_\tau(tw)[w,w]\dd t.
\]
For each fixed finite $k$, shrink the ball so that $D^j\widehat\Phi_\tau$ is bounded for $j\le\max\{2,k\}$. For $R(w)$ equal to the left-hand side, $R=O(\|w\|^2)$, $DR=O(\|w\|)$ and the higher derivatives are bounded. If $g_\varepsilon$ and all its derivatives through order $k$ are $O(\varepsilon)$, then each higher chain-rule term in $D^j(R\circ g_\varepsilon)$ is $O(\varepsilon^2)$: the single-$DR$ term contains an additional $O(\varepsilon)$, while every other term contains at least two derivatives of $g_\varepsilon$. This proves the stated $C^k$ remainder without requiring $f\in C^{k+1}$ or compactness of the infinite-dimensional stable ball.
\end{remark}

\begin{remark}[Finite regularity cannot be replaced by $C^\infty$]
The qualifier in Proposition~\ref{prop:analytic} is essential. Select an even, positive $p\in H^{s+4}(\T)\setminus C^\infty(\T)$ with $\int_\T p=1$ and $\|p-1\|_{H^{s+4}}$ arbitrarily small. Define
\[
 V=2\nu^2\frac{p''}{p}\in H^{s+2}_{\rm ev},
 \qquad Z_V=S_V(\tau)1,\qquad m=Z_Vp.
\]
Then $A_Vp=0$ and $S_V(\tau)p=p$. The mass of $m$ equals one, since $\int m=\langle S_V(\tau)1,p\rangle=\int S_V(\tau)p=1$; for $p$ close enough to $1$, $m$ is a positive density near $1$. Nevertheless $\Phi_\tau(V,m)=S_V(\tau)p=p\notin C^\infty$. Thus an unqualified claim of infinite regularization would be false within the same arbitrarily small potential ball.
\end{remark}

\section{A covering lemma with an infinite-dimensional contracting fibre}\label{app:cover}

This appendix gives the covering device used in Theorems~\ref{thm:symbolic} and \ref{thm:robust}. Its graph-transform components are standard (compare \cite{HPS,ZG}); the compact-core step is included because a closed Banach ball need not be compact.

Let $E$ be a Banach space, $I_1,\ldots,I_N\subset\R$ pairwise disjoint compact nondegenerate intervals, and $\mathcal N_i=I_i\times\overline B_\delta(E)$. For a map $T=(F,G)$ defined near $\mathcal N=\bigcup_i\mathcal N_i$, let
\[
 \begin{aligned}
 A&=\partial_xF, & B&=\partial_xG, & C&=D_yF, & D&=D_yG,\\
 \lambda_u&=\inf_{\mathcal N}|A|, & \beta&=\sup_{\mathcal N}\norm B, &
 \gamma&=\sup_{\mathcal N}\norm C, & \alpha&=\sup_{\mathcal N}\norm D.
 \end{aligned}
\]

\begin{lemma}[Compact shift factor and periodic closing]\label{lem:cover}
Suppose $T$ is $C^1$ on a neighbourhood of $\mathcal N$, its image $T(\mathcal N)$ is relatively compact, and for some $L>0$
\begin{equation}\label{eq:cones}
 \lambda=\lambda_u-\gamma L>1,\qquad
 \beta+\alpha L\le\lambda L,\qquad
 q=\alpha+\frac{(\beta+\alpha L)\gamma}{\lambda}<1,\qquad
 \sup_{\mathcal N}\norm G<\delta.
\end{equation}
Suppose also that for each $i$, $F(\cdot,y)$ has constant orientation on $I_i$ and, for every $\norm y\le\delta$, maps $I_i$ strictly across $\operatorname{conv}(\cup_jI_j)$ with uniform endpoint inequalities. Then there is a nonempty compact set $\Lambda_\infty\subset\mathcal N$ such that $T(\Lambda_\infty)=\Lambda_\infty$, and a continuous surjection from $\Lambda_\infty$ onto the full one-sided $N$-shift intertwining $T$ and the shift. Every periodic symbolic word of least period $n$ is realized by a point of least period $n$.
\end{lemma}
\begin{proof}
Let $\mathcal H_i$ be the complete metric space of Lipschitz maps $h:I_i\to\overline B_\delta(E)$ with Lipschitz constant at most $L$, with the uniform distance. Along the horizontal graph $(x,h(x))$, the scalar $F(x,h(x))$ is strictly monotone: by the fixed orientation of $\partial_xF$ and the Lipschitz bound on $h$, differences have absolute value at least $\lambda|x-x'|$. By the uniform endpoint covering it admits, for each target interval $I_j$, a unique inverse branch $\sigma_{h,j}:I_j\to I_i$ with $\mathrm{Lip}(\sigma_{h,j})\le\lambda^{-1}$. The transformed graph is
\[
  \mathcal T_{i\to j}h(x')=G(\sigma_{h,j}(x'),h(\sigma_{h,j}(x'))).
\]
The slope inequality in \eqref{eq:cones} gives $\mathrm{Lip}(\mathcal T_{i\to j}h)\le(\beta+\alpha L)/\lambda\le L$; stable trapping gives values in $\overline B_\delta(E)$. Comparing two inverse branches and using the bound for $C$ gives
\[
\norm{\sigma_{h,j}-\sigma_{\tilde h,j}}_\infty
\le\frac\gamma\lambda\norm{h-\tilde h}_\infty.
\]
The bound for $G$ then implies $\norm{\mathcal T_{i\to j}h-\mathcal T_{i\to j}\tilde h}_\infty\le q\norm{h-\tilde h}_\infty$.

To construct forward symbolic trajectories, start with an arbitrary horizontal graph over $I_{i_0}$. For every finite word $i_0\ldots i_n$, repeated inverse branches select a nonempty closed subset of the compact initial interval; these sets are nested when the word is extended. Their intersection is nonempty, and its orbit has the prescribed itinerary. Hence the closed survivor set
\[
 \Lambda=\{z\in\mathcal N:T^jz\in\mathcal N\ \text{for all }j\ge0\}
\]
is nonempty and realizes every one-sided sequence. Because the intervals $I_i$ are disjoint, its itinerary map $\pi$ to the full shift is continuous and intertwines $T$ with $\sigma$. The set $\Lambda$ is forward invariant but need not be compact, since it lies in a Banach ball. Define $K=\overline{T(\Lambda)}$. Relative compactness of $T(\mathcal N)$ makes $K$ compact. Since $T(\Lambda)\subset\Lambda$ and $\Lambda$ is closed, $K\subset\Lambda$ and $T(K)\subset K$. The nonempty nested intersection
\[
 \Lambda_\infty=\bigcap_{j\ge0}T^j(K)
\]
is compact and satisfies $T(\Lambda_\infty)=\Lambda_\infty$: the reverse inclusion follows by taking preimages in the nested nonempty compact sets $T^{j}(K)\cap T^{-1}\{z\}$ for $z\in\Lambda_\infty$. To see that every itinerary survives on $\Lambda_\infty$, fix a symbolic sequence $\boldsymbol{s}$ and write $K_j=T^j(K)$. For each $j$, prepend $j+1$ arbitrary symbols to $\boldsymbol{s}$, realize the longer sequence at a point of $\Lambda$, and advance that point $j+1$ times. The resulting point belongs to $K_j\cap\pi^{-1}\{\boldsymbol{s}\}$. These sets are nonempty, nested and compact, so their intersection is nonempty and is contained in $\Lambda_\infty\cap\pi^{-1}\{\boldsymbol{s}\}$. This compact-fibre argument is needed beyond separate realizability of finite prefixes. The factor inequality on this compact invariant restriction gives $h_{\rm top}(T|_{\Lambda_\infty})\ge\log N$ \cite{Walters}.

Finally, fix a periodic word $i_0\ldots i_{n-1}$ repeated indefinitely. The cyclic composition of the associated $n$ graph transforms is a contraction of $\mathcal H_{i_0}$ with factor at most $q^n$, so it has a fixed horizontal graph. Along the resulting cycle of graphs, the composite inverse scalar branch maps $I_{i_0}$ into itself and has Lipschitz constant at most $\lambda^{-n}$. Its fixed point gives a genuine $n$-periodic orbit. Its itinerary has least period $n$, so the point cannot have a smaller period. Since the entire orbit lies in $K$ and has an arbitrarily long prehistory along itself, it lies in $\Lambda_\infty$.
\end{proof}

\begin{remark}
A full-shift factor is weaker than a conjugacy. In particular, it does not establish that the invariant core is a Cantor set, that the compact core has a positive-measure basin, or that numerical survivor sets converge at a quantitative Hausdorff rate.
\end{remark}

\section{Comparison with the untuned long-horizon family}\label{app:long}

For the original horizon-independent construction set
\[
 V_\varepsilon[m]=-2\nu^2q_1^2\varepsilon f(c_1(m)/\varepsilon)e_1.
\]
This is not the interaction \eqref{eq:b}. For every fixed $\tau$, the first-order target of this family is $E_1a+(1-E_1)^2f(a)$, not $f(a)$. As the horizon grows, however, the normalized ground state of $A_V$ has first-order response $f$, and the following estimate makes this precise without a logarithmic horizon--amplitude constraint.

\begin{proposition}[Fixed long-horizon threshold]\label{prop:long}
Let $f\in C^1(U)$, $J\Subset U$, and $\nu>0$. In the scaled coordinates of \eqref{eq:coords}, for the horizon-independent $V_\varepsilon$, there exist $d,C>0$ and $\varepsilon_0>0$ such that for $0<\varepsilon<\varepsilon_0$ and all $\tau\ge1$,
\begin{equation}\label{eq:long}
 \norm{T_{\varepsilon,\tau}-(f,0)}_{C^1(\mathcal D)}
 \le C\{\varepsilon+(1+\tau)e^{-d\tau}\}.
\end{equation}
Thus any sufficiently large fixed horizon, followed by sufficiently small $\varepsilon$, preserves strictly covering expanding profiles.
\end{proposition}
\begin{proof}
We prove the parameter-uniform spectral estimate directly, rather than assuming it from an abstract eigenprojection theorem. Put $a_1=\nu q_1^2$ and $d=a_1/2$. Choose a sufficiently small potential ball in $H^{s+2}_{\rm ev}$ so that $\|V\|_\infty/(2\nu)\le a_1/8$. The min-max principle for the self-adjoint $A_V$ shows that its top eigenvalue $\ell(V)$ is simple and lies in $[-a_1/8,a_1/8]$, whereas all remaining even eigenvalues are at most $-7a_1/8$. In particular, the spectral gap is at least $3a_1/4$.

The augmented elliptic map
\[
 (V,p,\ell)\mapsto
\left(\nu p''-\frac{V}{2\nu}p-\ell p,\ \int_\T p-1\right)
\]
acts from $H^{s+2}_{\rm ev}\times H^{r+2}_{\rm ev}\times\R$ into $H^r_{\rm ev}\times\R$. Its derivative in $(p,\ell)$ at $(V,p,\ell)=(0,1,0)$ is $(h,\mu)\mapsto(\nu h''-\mu,\int h)$, which is invertible. The analytic implicit-function theorem supplies an analytic eigenpair $(p_V,\ell(V))$ with $\int p_V=1$. After reducing the potential ball, $p_V$ is positive and bounded near $1$ in $H^{r+2}$, and its first Fr\'echet derivative is uniformly bounded. This argument does not use compactness of the potential ball.

Let $\Pi_Vh=p_V\langle p_V,h\rangle_{L^2}/\|p_V\|_2^2$, $Q_V=I-\Pi_V$, $B_V=A_V-\ell(V)I$, and $E_V(t)=e^{tB_V}Q_V$. By self-adjointness,
\[
\|E_V(t)\|_{L^2\to L^2}\le e^{-3a_1t/4}.
\]
We need this uniformly on $H^r$, not merely on $L^2$. Bounded perturbation gives uniform $H^r$ semigroup norms on bounded time intervals. Uniform smoothing from $L^2$ into $H^r$ at any fixed positive time follows by successive heat/Duhamel estimates, using the $H^{s+2}$ multiplier bound. Splitting the semigroup into two fixed smoothing pieces and an $L^2$ piece for $t\ge1$ therefore gives
\begin{equation}\label{eq:decayE}
\|E_V(t)\|_{H^r\to H^r}\le Ce^{-dt},\qquad t\ge0.
\end{equation}
The differentiated estimate must retain the moving spectral projection. Since $Q_V$ commutes with $B_V$,
\[
 D_VE_V(t)[H]=D_VQ_V[H]E_V(t)+E_V(t)D_VQ_V[H]
+\int_0^t E_V(t-u)D_VB_V[H]E_V(u)\,\dd u.
\]
Here $D_VB_V[H]=-M_H/(2\nu)-D\ell(V)[H]I$ and the derivatives of $Q_V$ and $\ell(V)$ have operator bounds $C\|H\|_{H^{s+2}}$. Combining with \eqref{eq:decayE} yields
\begin{equation}\label{eq:ddecayE}
\|D_VE_V(t)[H]\|_{H^r\to H^r}
\le C(1+t)e^{-dt}\|H\|_{H^{s+2}}.
\end{equation}
The projection-derivative terms are indispensable: differentiating only the restricted semigroup would not justify uniform dependence on $V$.

Define
\[
 b_V=\Pi_V1=\frac{p_V}{\|p_V\|_2^2},\qquad
 z_{V,t}=e^{-t\ell(V)}S_V(t)1=b_V+E_V(t)1.
\]
Since $b_0=1$ and $Q_V1=1-b_V=O(\|V\|_{H^{s+2}})$, both $z_{V,t}$ and $b_V$ stay uniformly close to $1$ in $H^r$ for all $t\ge0$. Their reciprocals and first parameter derivatives are uniformly bounded in the Banach algebra $H^r$ on a smaller fixed ball. Set
\[
\mathcal L_{V,t}=(\Pi_V+E_V(t))M_{1/z_{V,t}},
\qquad \mathsf P_V=\Pi_VM_{1/b_V}.
\]
For $h\in H^r$, the normalization $\int p_V=1$ gives $\mathsf P_Vh=p_V\int_\T h$. Hence
\[
\mathcal R_{V,t}:=\mathcal L_{V,t}-\mathsf P_V
=\Pi_VM_{1/z_{V,t}-1/b_V}+E_V(t)M_{1/z_{V,t}},
\qquad \frac1z-\frac1b=-\frac{E_V(t)1}{zb}.
\]
The algebra estimates together with \eqref{eq:decayE}--\eqref{eq:ddecayE} give the uniform, differentiated quotient bounds
\begin{equation}\label{eq:rem-long}
\|\mathcal R_{V,t}\|_{\mathcal L(H^r)}\le Ce^{-dt},\qquad
\|D_V\mathcal R_{V,t}[H]\|_{\mathcal L(H^r)}
\le C(1+t)e^{-dt}\|H\|_{H^{s+2}},\quad t\ge1.
\end{equation}

For the untuned feedback of the proposition, put $m=1+\varepsilon(ae_1+y)$ and $V=-2\nu^2q_1^2\varepsilon f(a)e_1$. The exact identity $\mathcal R_{0,t}1=0$ is important. It gives
\[
\mathcal R_{V,t}m=\mathcal R_{V,t}(m-1)
 +\int_0^1D_V\mathcal R_{\theta V,t}[V]1\,\dd\theta
=O\bigl(\varepsilon(1+t)e^{-dt}\bigr).
\]
For a scaled direction $\xi$, one differentiates directly:
\[
D_\xi(\mathcal R_{V,t}m)
=D_V\mathcal R_{V,t}[D_\xi V]m+\mathcal R_{V,t}D_\xi m.
\]
Here $D_\xi V$ and $D_\xi m$ are $O(\varepsilon)$, so the same uniform bound follows in $C^1$ without requiring $f\in C^2$. Finally, from the normalized eigenfunction equation,
\[
D_Vp_0[e_1]=-\frac{e_1}{2\nu^2q_1^2},\qquad
p_{-2\nu^2q_1^2\varepsilon f(a)e_1}
=1+\varepsilon f(a)e_1+O_{C^1}(\varepsilon^2).
\]
Combining $\Phi_t(V,m)=p_V+\mathcal R_{V,t}m$ with \eqref{eq:rem-long} and the preceding bounds, projecting onto $e_1$ and its complement, and dividing by $\varepsilon$ proves \eqref{eq:long}. The constants can be chosen uniformly for all $t=\tau\ge1$ on the fixed scaled tube; no claim is made uniformly as $\tau\downarrow0$.
\end{proof}

\section{Optimal-control interpretation and limitations}\label{app:control}
For a fixed positive input density $m$ with frozen running potential $V[m]$, consider the controlled diffusion on the torus
\[
\dd X_t=\alpha_t\,\dd t+\sqrt{2\nu}\,\dd B_t,\qquad X_0\sim m,
\]
with running cost $\tfrac12|\alpha_t|^2+V[m](X_t)$ and zero terminal cost. If $u$ solves \eqref{eq:PDE}, It\^o's formula and completion of the square imply
\[
\mathbb E\int_0^\tau\left(\tfrac12|\alpha_t|^2+V[m](X_t)\right)\dd t
=\mathbb E\,u(0,X_0)
+\frac12\mathbb E\int_0^\tau|\alpha_t+u_z(t,X_t)|^2\,\dd t.
\]
The stochastic integral has zero mean for square-integrable admissible controls after standard localization; on the compact torus the derivatives of $u$ supplied by the parabolic regularity above are bounded. Equality in the value inequality is achieved at $\alpha_t=-u_z(t,X_t)$, and its law solves the forward equation in \eqref{eq:PDE}. The Feynman--Kac representation for $S_V(t)$ has a weight bounded by $e^{\pm t\|V\|_\infty/(2\nu)}$, so no extra integrability condition is hidden. This verifies the control interpretation of the frozen response; it is not a Nash-consistency condition for a time-evolving population and does not establish a finite-player approximation.

There is a concrete functional obstruction to a direct finite-player interpretation. An empirical probability measure is atomic and does not belong to $H^r(\T)$ for $r>1/2$, whereas all the differentiability arguments above concern strictly positive Sobolev densities. A particle model would need a topology admitting empirical measures, perhaps after smoothing, together with one-episode approximation estimates. Even if such an approximation incurred error $\xi_N$ at each episode and the law response had Lipschitz constant $L>1$, the repeated-error recurrence $\Delta_{n+1}\le L\Delta_n+\xi_N$ would yield only
\[
\Delta_n\le L^n\Delta_0+\xi_N\frac{L^n-1}{L-1}.
\]
Thus a small one-step error alone does not control arbitrarily many expanding episodes. Invariant-law stochastic stability or a quantitative many-particle limit would require a separate analysis.

\section*{Declaration of competing interest}
All authors declare no competing interest.

\section*{Acknowledgments}
This work was supported by the National Natural Science Foundation of
China (Grant No. 12461054), and the Guizhou Province High-Level Talent
Cultivation Support Program (Grant No. KJLYRC-[2026]014). The authors
acknowledge this support.

\section*{Data Availability}
Data will be made available on request.

\end{document}